%  AMSTeX
%  Dec. 1, 1999
\documentstyle{amsppt}
\magnification=1200
\hsize=6.5truein
\vsize=8.9truein
\topmatter
\title Algebraic Gauss-Manin Systems \\
and Brieskorn Modules
\endtitle
\author   Alexandru Dimca and Morihiko Saito
\endauthor
\abstract
We study the algebraic Gauss-Manin system and the algebraic 
Brieskorn module associated to a polynomial mapping with isolated
singularities.
Since the algebraic Gauss-Manin system does not contain any
information on the cohomology of singular fibers, we first construct a 
non quasi-coherent sheaf which gives the cohomology of every fiber.
Then we study the algebraic Brieskorn module, and show that its 
position in the the algebraic Gauss-Manin system is determined by a 
natural map to quotients of local analytic Gauss-Manin systems, and 
its pole part by the vanishing cycles at infinity, comparing it with the Deligne extension.
This implies for example a formula for the determinant of periods.
In the two-dimensional case we can describe the global structure of 
the algebraic Gauss-Manin system rather explicitly.
\endabstract
\endtopmatter
\tolerance=1000
\baselineskip=12pt 

\document
\NoBlackBoxes

\def\scirc{\raise.2ex\hbox{${\scriptstyle\circ}$}}
\def\ssbull{\raise.2ex\hbox{${\scriptscriptstyle\bullet}$}}
\def\cbull{\raise.5ex\hbox{${\scriptscriptstyle\bullet}$}}

\def\fin{\text{\rm fin}}

\def\ord{\hbox{{\rm ord}}}
\def\res{\hbox{{\rm res}}}
\def\Kernel{\hbox{{\rm Kernel}}}
\def\Tr{\hbox{{\rm Tr}}}
\def\Cone{\hbox{{\rm Cone}}}
\def\Perv{\hbox{{\rm Perv}}}

\def\Gr{\text{{\rm Gr}}}
\def\Im{\hbox{{\rm Im}}}

\def\Coker{\hbox{{\rm Coker}}}
\def\Ker{\hbox{{\rm Ker}}}

\def\GAGA{\hbox{{\rm GAGA}}}
\def\DR{\hbox{{\rm DR}}}
\def\IC{\hbox{{\rm IC}}}

\def\supp{\hbox{{\rm supp}}\,}
\def\Sing{\text{{\rm Sing}}\,}

\def\an{\text{\rm an}}
\def\rank{\text{\rm rank}}

\def\tor{\text{\rm tor}}

\def\simto{\buildrel\sim\over\to}

\def\SameAuthor{\vrule height3pt depth-2.5pt width1cm}

\bigskip
\centerline{\bf Introduction}

\bigskip\noindent
Let
$ f : X \rightarrow  S $ be a morphism of smooth complex algebraic 
varieties.
The Gauss-Manin systems
$ {\Cal G}_{f}^{i} $ of
$ f $ are defined to be the (cohomological) direct images of
$ {\Cal O}_{X} $ as algebraic left
$ {\Cal D} $-Modules.
They are regular holonomic
$ {\Cal D}_{S} $-Modules, and correspond by the de Rham functor
$ \DR $ to the (perverse) higher direct images of the constant sheaf.
In particular, they give the cohomology with compact support of each 
fiber,
but not the cohomology of singular fibers unless $ f $ is proper, 
because the stalk of the higher direct image does not coincide with the
cohomology of the fiber due to the vanishing cycles at infinity.
Furthermore the cohomologies of the fibers do not form a constructible 
sheaf in a natural way (see (2.6)), 
and it is not easy to construct a quasi-coherent sheaf on 
$ S $ which is generically coherent, 
and contains the information on the cohomology of every fiber.

Assume for simplicity
$ X = {\Bbb A}^{n}, S = {\Bbb A}^{1} $ with
$ n \ge  2 $.
Then it is easy to show that
$$
\DR_{S}({\Cal G}_{f}^{i}) = (R^{n+i-1}f_{*}{\Bbb C}_{X})[1]
$$
in the category of perverse sheaves
$  \Perv(S,{\Bbb C}) $ (see [3]).
Let
$ {\Cal G}_{f} = {\Cal G}_{f}^{0} $,
$ G_{f} = \Gamma (S, {\Cal G}_{f}) $, and
$ L = R^{n-1}f_{*}{\Bbb C}_{X} $.
Note that
$$
{\Cal G}_{f}^{i} = 0\quad\text{for }i \ne 1-n, 0,
$$
if
$ f $ has at most isolated singularities including at infinity [31]
(or, more generally, if the support of the vanishing cycles including at infinity is discrete).  See (1.3).

For a variety
$ Y $ of pure dimension
$ m $ in general,
we have the trace morphism
$ \Tr : {H}_{c}^{2m}(Y,{\Bbb C}) \rightarrow  {\Bbb C} $,
and
$ {\tilde{H}}_{c}^{i}(Y,{\Bbb C}) $ is defined as
$ \Ker(\Tr : {H}_{c}^{2m}(Y,{\Bbb C}) \rightarrow  {\Bbb C}) $ for
$ i = 2m $,
and
$ {H}_{c}^{i}(Y,{\Bbb C}) $ otherwise.

For
$ c \in  {\Bbb C} $,
let
$ X_{c} = f^{-1}(c) $.
If
$ f $ has at most isolated singularities including at infinity, then
$ {\tilde{H}}_{c}^{i}(X_{c},{\Bbb C}) = 0 $ for
$ i \ne  n - 1, n $,
and we can easily show canonical isomorphisms
$$
\Ker(t-c|G_{f}) = {\tilde{H}}_{c}^{n}(X_{c},{\Bbb C})^{*},\quad 
\Coker(t-c|G_{f}) = {H}_{c}^{n-1}(X_{c},{\Bbb C})^{*},
$$
where
$ * $ denotes the dual vector space.
See (1.2--3).

From now on, we assume that
$ f $ has at most isolated singularities.
Let
$ \Omega ^{i} = \Gamma (X,{\Omega }_{X}^{i}) $.
Then there exists a natural morphism
$ \Omega ^{n} \rightarrow  G_{f} $,
and its kernel is
$ df\wedge d\Omega ^{n-2} $.
So we have a
$ {\Bbb C}[t] $-submodule
$$
{G}_{f}^{(0)} = \Omega ^{n}/df\wedge d\Omega ^{n-2} \subset  G_{f},
$$
which is called the {\it algebraic Brieskorn module} of
$ f $.
Let
$ {\Cal G}_{f}^{(0)} $ denote the quasi-coherent sheaf corresponding to
$ {G}_{f}^{(0)} $,
and
$ {\Cal G}_{f}^{(0),\an} $ the associated analytic sheaf.

For
$ x \in  X $,
let
$ {\Cal G}_{f,x}^{(0)} $ denote the local analytic Brieskorn lattice
$ {\Omega }_{{X}^{\an},x}^{n}/df\wedge d{\Omega }_{{X}^{\an},x}^{n-2} $ 
(see [7]), and 
$ {\Cal G}_{f,x} $ the local analytic Gauss-Manin system [21]
which is the localization of
$ {\Cal G}_{f,x}^{(0)} $ by
$ {\partial }_{t}^{-1} $.
(They vanish unless
$ x \in   \Sing f $.)
Then for
$ x \in  X_{c} $,
we have the restriction morphism
$ {\Cal G}_{f,c}^{(0),\an} \rightarrow  {\Cal G}_{f,x}^{(0)} $.
We define
$$
{\Cal L}_{f,c} = \Ker({\Cal G}_{f,c}^{(0),\an} \rightarrow  
\bigoplus _{f(x) =c} {\Cal G}_{f,x}^{(0)}).
$$

\proclaim{\bf 0.1.~Theorem}
If
$ f $ has at most isolated singularities including at infinity,
we have natural isomorphisms 
$$
\Ker(t-c|{\Cal L}_{f,c}) = \tilde{H}^{n-2}(X_{c},{\Bbb C}),\quad 
\Coker(t-c|{\Cal L}_{f,c}) = H^{n-1}(X_{c},{\Bbb C}),
$$
and 
$ \tilde{H}^{i}(X_{c},{\Bbb C}) = 0 $ for  
$ i \ne  n - 2, n - 1$.
See (2.4).
\endproclaim

But the
$ {\Cal L}_{f,c} $ do not form a quasi-coherent sheaf on
$ S^{\an},  $ and
$ {\Cal G}_{f,c}^{(0),\an} $ in the definition of
$ {\Cal L}_{f,c} $ cannot be replaced by
$ {\Cal G}_{f,c}^{(0)} $ or
$ {G}_{f}^{(0)} $.
(Indeed, the image of the natural morphism 
$ G_{f}^{(0)}\to {\Cal G}_{f,x}^{(0)} $ is not a 
$ {\Bbb C}[t] $-module of rank
$ \mu_{x} $ in general, where  
$ \mu _{x} $ denote the Milnor number of
$ f $ at
$ x \in   \Sing f $.)
From (0.1) we get

\bigskip

\noindent
{\bf 0.2.~Corollary.}
With the assumption of (0.1) we have
$$
\aligned
\dim \Ker(t-c|{G}_{f}^{(0)}) &= \dim \tilde{H}^{n-2}(X_{c},{\Bbb C}),
\\
\dim \Coker(t-c|{G}_{f}^{(0)}) &= \dim H^{n-1}(X_{c},{\Bbb C}) + 
\sum _{f(x) = c} \mu _{x}.
\endaligned
$$
See also (3.6) for the relation with other invariants.

In the local analytic case, the Brieskorn lattice
$ {\Cal G}_{f,x}^{(0)} $ gives the Hodge numbers of the local analytic 
Milnor fiber [32] by [35] (see also [22], [27], [29], etc.), and contains the 
information on the analytic structure of
$ (f,x) $ in the local moduli space (see for example [28]).
However, concerning the global algebraic structure of
$ f $,
the algebraic Brieskorn module
$ {G}_{f}^{(0)} $ does not contain any more than the information on the local 
analytic structure of
$ (f,x) $ at
$ x \in   \Sing  f $ (although it gives the sums of the Hodge numbers
of local Milnor fibers, see (3.5)),
because
$ {G}_{f}^{(0)} $ is determined only by
$ G_{f} $ together with the composition of natural morphisms
$ {\Cal G}_{f,c}^{\an} \to {\Cal G}_{f,x} \to 
{\Cal G}_{f,x}/{\Cal G}_{f,x}^{(0)} $
for
$ x \in \Sing f $ due to the following :

\proclaim{\bf 0.3.~Theorem}
$$
{G}_{f}^{(0)} = \Ker\bigl(G_{f} \rightarrow  \bigoplus _{x\in \Sing f} 
{\Cal G}_{f,x}/{\Cal G}_{f,x}^{(0)}\bigr).
$$
\endproclaim

\noindent
(This is a part of Theorem (0.5) below.)
Note that
$ G_{f} $ is determined by the constructible sheaf
$ L \,(:= R^{n-1}f_{*}{\Bbb C}_{X}) $ by using the Riemann-Hilbert correspondence (see for example [4]),
and the morphism
$  {\Cal G}_{f,c}^{\an} \to {\Cal G}_{f,x} $ for
$ x \in X_{c} $ is determined also topologically by
using  the restriction to the local Milnor fibration.
It is relatively easy to describe
$ {\Cal G}_{f}^{\an} $ at least locally,
using the local classification of regular holonomic
$ {\Cal D} $-module of one variable [5], [6] (see also [26, 1.3]).

From (0.3) we also get

\proclaim{\bf 0.4.~Corollary}
$$
{\Cal L}_{f,c} = \Ker\bigl({\Cal G}_{f,c}^{\an} \rightarrow  
\bigoplus _{f(x)=c} {\Cal G}_{f,x}\bigr).
$$
\endproclaim

In particular,
$ {\Cal L}_{f,c} $ is actually a
$ {D}_{S,c}^{\an} $-module (on which the action of
$ \partial _{t} $ is surjective).
The corresponding constructible sheaf defined on a neighborhood 
$ S' $  of 
$ c $ is given by 
$ R^{n-1}f_{*}\tilde{j}_{!}{\Bbb C}_{X\setminus B}|_{S'} $ where 
$ B = \bigcup _{x\in \Sing f} B_{x} $ with 
$ B_{x} $ a sufficiently small ball with center
$ x $,  and 
$ \tilde{j}  : X \setminus B \to X $ denotes the inclusion.

In general,
$ {G}_{f}^{(0)} $ is not finitely generated over
$ {\Bbb C}[t] $.
As a corollary of Sabbah's results [23],
$ {G}_{f}^{(0)} $ is a finite
$ {\Bbb C}[t] $-module if and only if
$ f $ has a certain good property at infinity (i.e.
if
$ f $ is cohomologically tame in his sense).
See (3.3).
Let
$ U $ be a dense open subvariety of
$ S $ such that
$ L|_{U} $ is a local system (i.e.
 $ {\Cal G}_{f}|_{U} $ is a locally free
$ {\Cal O}_{U} $-Module of finite rank).
Put
$ \Delta = S \setminus U $.
Let
$ {\Cal G}_{f}^{>0} $ be the Deligne extension of
$ {\Cal G}_{f}|_{U} $ such that the eigenvalues of the residues of the 
connection are contained in
$ (-1, 0] $.
See [9].
(The shift of the index comes from the normalization of the exponents
in [25].)
Put
$ {G}_{f}^{>0} = \Gamma (S, {\Cal G}_{f}^{>0}) $.
Then
$ {G}_{f}^{>0} $ is always a free
$ {\Bbb C}[t] $-module of finite type.
We will see as a corollary of (0.5) below that
$ {G}_{f}^{(0)} $ is finite over
$ {\Bbb C}[t]  $ if and only if it is contained in
$ {G}_{f}^{>0} $.
Let
$$
{G}_{f}^{(0),>0} = {G}_{f}^{(0)} \cap  {G}_{f}^{>0}.
$$
This is called the {\it Brieskorn-Deligne lattice}.
It is a free
$ {\Bbb C}[t] $-module of finite rank, and generates
$ G_{f} $ over
$ \Gamma (S,{\Cal D}_{S}) $,
because
$ {\Cal G}_{f} $ has no nontrivial quotient with discrete support.
See (1.3).
For
$ x \in   \Sing  f $,
we can define the Deligne extension
$ {\Cal G}_{f,x}^{>0} $ (which is contained in
$ {\Cal G}_{f,x}) $ similarly.
Then

\proclaim{\bf 0.5.~Theorem}
We have the following commutative diagram of exact sequences of
$ {\Bbb C}[t] $-modules.
See (3.2).
\endproclaim

$$
\CD
@. 0 @. 0 @. 0 \\
@. @VVV @VVV @VVV \\
0 @>>> {G}_{f}^{(0),>0} @>>> {G}_{f}^{>0} @>>>
\bigoplus _{x\in \Sing f} {\Cal G}_{f,x}^{>0}/{\Cal G}_{f,x}^{(0)} @>>> 0 
\\
@. @VVV @VVV @VVV \\
0 @>>> {G}_{f}^{(0)} @>>> {G}_{f} @>>>
\bigoplus _{x\in \Sing f} {\Cal G}_{f,x}/{\Cal G}_{f,x}^{(0)} @>>> 0 \\
@. @VVV @VVV @VVV \\
0 @>>> \bigoplus _{c\in \Delta } P_{c}(E'_{c},T) @>>> 
\bigoplus _{c\in \Delta } P_{c}(E_{c},T) @>>>
\bigoplus _{x\in \Sing f} P_{f(x)}(E''_{x},T) @>>> 0 \\
@. @VVV @VVV @VVV \\
@. 0 @. 0 @. 0 \\
\endCD
$$
See (3.1) for
$ P_{c}(E'_{c},T),  $ etc., where
$ c $ and 
$ t - c $ are denoted by 
$ s $  and
$ t_{s} $ respectively.
Note that
$ P_{c}(E'_{c},T) $ is isomorphic to
$ E'_{c}\otimes _{{\Bbb C}}{\Bbb C}[{1 \over t-c}]{1 \over t-c} $ as a
$ {\Bbb C}[t] $-module (where
$ {\Bbb C}[{1 \over t-c}]{1 \over t-c} := {\Bbb C}[t,{1 \over t-c}]/{\Bbb C}[t]) $,
and is called the pole part of
$ {G}_{f}^{(0)} $ at
$ c $.
So every vanishing cycle at infinity contributes to the 
pole part of
$ {G}_{f}^{(0)} $, since
$ E'_{c} $ is the space of vanishing cycles at infinity.
For the proof of (0.5), it is not necessary to use the theory of the direct 
image of the filtration 
$ V $ of Kashiwara [18] and Malgrange [20] as in [23], [24], [27],
because it is actually enough to use a much easier theory of the filtration
$ V $ in the one variable case which is closely related to
the Deligne extension [9] and also to Varchenko's theory [35], [36].

Theorem (0.5) describes the structure of
$ {G}_{f}^{(0)} $ by comparing it with the Deligne extension
$ {G}_{f}^{>0} $.
The first row means that
$ {G}_{f}^{(0),>0} $ coincides with
$ {G}_{f}^{>0} $ up to torsion, because we have for
$ x\in \Sing X_{c} $
$$
{\Cal G}_{f,x}^{>0}/{\Cal G}_{f,x}^{(0)} = \sum _{1\le j\le \mu _{x}} 
{\Bbb C}[t]/(t - c)^{k_{j} }{\Bbb C}[t]
$$
as
$ {\Bbb C}[t] $-modules.
Here
$ k_{j} $ are nonnegative integers such that
$ 0 < \alpha _{j} - k_{j} \le  1 $ with
$ \alpha _{j}  \,(1   \le  j \le  \mu _{x}) $ the exponents of the 
Brieskorn lattice.
See [26].
Let
$ {\mu }_{x}^{1} $ be the dimension of the unipotent monodromy part of 
the vanishing cohomology of
$ f $ at
$ x $ (i.e.
the number of the exponents which are integers.)
Then from the first row of (0.5) together with the symmetry of the 
exponents, we get

\bigskip
\noindent
{\bf 0.6.~Corollary.} 
$$
\dim {\Cal G}_{f,c}^{>0}/{\Cal G}_{f,c}^{(0),>0} = 
\sum _{f(x)=c} {1 \over 2}((n-1)\mu _{x} - {\mu }_{x}^{1}).
$$

In particular,
$ {G}_{f}^{(0),>0} = {G}_{f}^{>0} $ if and only if
$ n = 2 $ and every singular point of
$ f $ on
$ X $ is an ordinary double point (i.e., a node), because the minimal 
exponent has multiplicity one.

Let
$ {G}_{f}^{(-1)} = df\wedge \Omega ^{n-1}/df\wedge d\Omega ^{n-2} 
\, (= \partial_{t}^{-1}{G}_{f}^{(0)})
\subset  G_{f}^{(0)} $ (similarly for
$ {\Cal G}_{f,x}^{(-1)} $), and
$ {G}_{f}^{(-1),>0} = {G}_{f}^{(-1)} \cap {G}_{f}^{>0} $.
Then (0.5) holds with
$ {G}_{f}^{(0)} $,
$ {G}_{f}^{(0),>0} $,
$ {\Cal G}_{f,x}^{(0)} $ replaced by
$ {G}_{f}^{(-1)} $,
$ {G}_{f}^{(-1),>0} $,
$ {\Cal G}_{f,x}^{(-1)} $.

Let
$ \{\omega _{i}\} $,
$ \{df\wedge\eta _{i}\} $ be
$ {\Bbb C}[t] $-bases of
$ {G}_{f}^{(0),>0} $,
$ {G}_{f}^{(-1),>0} $,
and
$ \{\gamma _{j}(c)\} $ a basis of horizontal (multivalued) sections of
$ \prod _{c\in U} H_{n-1}(X_{c},{\Bbb Z})/\text{torsion} $.
Then the squares of the determinants of the period matrices
$$
\det\biggl(\int_{{\gamma }_{j}(t)}\res{{\omega }_{i} \over f-t}\biggr)^{2},
\quad  
\det\biggl(\int_{{\gamma }_{j}(t)}{\eta }_{i} \biggr)^{2}
\quad  \text{for }t \in  U
$$
are holomorphic functions on
$ U $,
and are independent of the choice of the bases up to constant multiples.
Here
$ \res{{\omega }_{i} \over f-t} $ denotes the Poincar\'e residue.
(See for example, [7,1.5]).
For
$ c \in  \Delta  := S \setminus  U $,
let
$ {\nu }_{c}^{\ne 1} $ be the dimension of the nonunipotent monodromy 
part of
$ E'_{c} $,
and put
$$
m_{c} = - {\nu }_{c}^{\ne 1} + \sum _{f(x) =c} (n - 2)\mu _{x},
\quad
m'_{c} = - {\nu }_{c}^{\ne 1} + \sum _{f(x) =c} n\mu _{x}.
$$
Then from (0.5) we can deduce (see (3.4)) :

\bigskip
\noindent
{\bf 0.7.~Corollary.}  
$$
\aligned
\det\biggl(\int_{{\gamma }_{j}(t)}^{}\res{{\omega }_{i} \over f-t}\biggr)^{2}
& =  \text{const }   \prod _{c\in \Delta } (t - c)^{m_{c}},
\\
\det\biggl(\int_{{\gamma }_{j}(t)}{\eta }_{i} \biggr)^{2}
& =  \text{const }   \prod _{c\in \Delta } (t - c)^{m'_{c}}.
\endaligned
$$

\noindent
This generalizes [16, 2.2] where
$ n = 2 $ and
$ f $ is a semiweighted homogeneous polynomial in the sense of loc.~cit.
(in particular,
$ f $ is tame and
$ \nu  = 0) $.
The corresponding assertion for the local analytic Brieskorn lattice is 
well-known.
See [36].

Since the algebraic Brieskorn module is determined by the algebraic
Gauss-Manin system
$ G_{f} $ together with the map to the quotients of the local analytic
Gauss-Manin systems, we consider the problem of determining the 
global structure of the Gauss-Manin system.
We treat rather the perverse sheaf
$$
L[1] := (R^{n-1}f_{*}{\Bbb C}_{X})[1]
$$
corresponding to
$ {\Cal G}_{f} $ by the Riemann-Hilbert correspondence,
because it is easier to handle.
In the case
$ n = 2 $,
we can describe it rather explicitly by using the relative version of
Deligne's weight spectral sequence [11] (see [24], [25]).

Let
$ \bar{f} : \bar{X} \rightarrow  S $ be a relative 
compactification of
$ f $ such that
$ \bar{X} $ is smooth and
$ \bar{L} := R^{1}\bar{f}_{*}{\Bbb C}_{\bar{X}}|_{U} $ is a 
local system (shrinking
$ U $ if necessary).
Let
$ g $ be the genus of the generic fiber of
$ \bar{f} $ such that
$ \rank \,\bar{L} = 2g $.
Let
$ j : U \rightarrow X $ be the inclusion morphism.
Then
$ (j_{*}\bar{L})[1] $ is the intersection complex
$ \IC_{S}\bar{L} $ (i.e. the intermediate direct image of
$ \bar{L}[1] $ in the sense of [3]).
See also [37].
Let
$ h $ be the number of horizontal irreducible components of
$ \bar{X} \setminus  X $, and
$ n_{c} $ the number of irreducible components of
$ X_{c} $.
Let
$ W $ be the weight filtration on
$ L[1] $ in the theory of mixed Hodge Modules [24] [25]
(this can also be obtained by using a mod
$ p $ reduction as in [3]).
Then, using the weight spectral sequence, we can easily show the following :

\proclaim{\bf 0.8.~Proposition}
With the above notation and assumption,
$ {\Gr}_{k}^{W}(L[1]) = 0 $ for
$ k \ne  2, 3 $,
and
$$
\aligned
{\Gr}_{2}^{W}(L[1]) &= \biggl(\bigoplus _{c\in \Delta } 
\biggl( {\overset{n_{c}-1}\to \bigoplus} \,
{\Bbb C}_{\{c\}}\biggr)\biggr)
 \bigoplus  (j_{*}\bar{L})[1],
\\
{\Gr}_{3}^{W}(L[1]) &= \biggl(\overset{h-1}\to \bigoplus \,
{\Bbb C}_{S}[1]\biggr) \bigoplus  
\biggl(\bigoplus _{i} (j_{*}L'_{i})[1]\biggr),
\endaligned
$$
where the
$ L'_{i} $ are non constant irreducible local systems on
$ U $.
\endproclaim

(The proof is more or less standard, and is left to the reader.
Indeed, the vanishing of
$ {\Gr}_{4}^{W} (L[1]) $ follows from (1.3), and the multiplicity 
$ n_{c}-1 $ of
$ {\Bbb C}_{\{c\}} $ is determined by using for example (1.2).)
By the long exact sequence associated with the filtration
$ W $ on
$ L[1] $, we get a refinement of Kaliman's inequality [17]:

\proclaim{\bf 0.9.~Corollary}
$$
h - 1 = \sum _{c\in \Delta } (n_{c} - 1) + \dim H^{1}(S^{\an}, 
j_{*}\bar{L}).
$$
\endproclaim

Passing to the corresponding regular holonomic
$ {\Cal D} $-modules, we get submodules
$ M' $,
$ M'' $ of
$ {\Cal G}_{f} $ with a short exact sequence
$$
0 \rightarrow  M' \bigoplus  M'' \rightarrow  {\Cal G}_{f} \rightarrow  
\tilde{M}\rightarrow  0,
$$
such that
$ M' = \bigoplus _{c} \bigl(\overset{n_{c}-1}\to\bigoplus\,
 {\Cal D}_{S}/{\Cal D}_{S}(t-c)\partial _{t} \bigr) $,
 $ \tilde{M} = \bigoplus _{i} \tilde{M}_{i} $ with
$ \DR_{S}(\tilde{M}_{i}) = (j_{*}L'_{i})[1] $.
As to
$ M'' $,
we have a short exact sequence
$$
0 \rightarrow  M_{S}(\bar{L}) \rightarrow  M'' \rightarrow  
\overset{h''}\to\bigoplus \,
{\Cal O}_{S} \rightarrow  0,
$$
such that
$ \DR_{S}(M_{S}(\bar{L})) = (j_{*}\bar{L})[1] $,
where
$ h'' = \dim H^{1}(S^{\an}, j_{*}\bar{L}) $.
So the description of
$ {\Cal G}_{f} $ is essentially reduced to the problem of extension.
For example, if
$ r = 0 $ (i.e.
 $ f : Y_{j} \rightarrow  S $ is bijective for any irreducible component
$ Y_{j} $ of
$ \bar{X} \setminus  X $ such that
$ \bar{f}(Y_{j}) = S) $,
then we have
$ {\Cal G}_{f} = M' \oplus  M'' $.
Note that
$ r = 0 $ by [15] if
$ n = 2 $ and
$ f $ is cohomologically tame [23] (see also [2]).

Note that (0.8) implies a nontrivial assertion on the local system
$ L|_{U} $.
In terms of the corresponding
$ {\Cal D} $-modules, let
$ {}_{\tor}{\Cal G}_{f} $ be the maximal
$ {\Cal D}_{S} $-submodule of
$ {\Cal G}_{f} $ whose support is contained in
$ \Delta  $.
Then we have
$$
{}_{\tor}{\Cal G}_{f} = \bigoplus _{c} \biggl(\overset{n_{c}-1}
\to\bigoplus \,
{\Cal D}_{S}/{\Cal D}_{S}(t-c) \biggr) \subset  M',
\quad \DR_{S}({}_{\tor}{\Cal G}_{f}) = 
\bigoplus _{c} \biggl( \overset{n_{c}-1}\to\bigoplus\, {\Bbb C}_{\{c\}} 
\biggr),
$$
and
$ M'/{}_{\tor}{\Cal G}_{f} \,(= \overset{h'}\to\bigoplus\, {\Cal O}_{S}) $ 
is a direct factor of
$ {\Cal G}_{f}/_{\tor}{\Cal G}_{f} $,where
$ h' = \sum _{c} (n_{c} - 1) $.
(This follows from the semisimplicity of
$ {\Gr}_{3}^{W}{\Cal G}_{f} $.)
It implies a result of Bailly-Maitre [2] that the maximal constant 
subsheaf of
$ L|_{U} $ has rank
$ h' $ (using [14]), and is a direct factor of
$ L|_{U} $.

In Sect. 1, we introduce the algebraic Gauss-Manin system, and study the action of 
$ t - c $ on it.
In Sect. 2, we define a filtration on the algebraic Gauss-Manin system
which gives the algebraic Brieskorn module, and prove (0.1--2) in a 
more general situation.
Then Theorem (0.5) and its corollaries are proved in Sect. 3.

\bigskip

\noindent
{\it Convention.} In this paper, algebraic variety means a separated 
scheme of finite type over
$ {\Bbb C} $.
For a variety
$ X $ and a morphism
$ f $,
we denote
$ {\Bbb C}_{X^{\an}} $ by
$ {\Bbb C}_{X} $,
and
$ f^{\an} $ by
$ f $ to simplify the notation, where
$ X^{\an} $ is the associated analytic space.
Similarly,
$ H^{i}(X,{\Bbb C}) $ and
$ {H}_{c}^{i}(X,{\Bbb C}) $ denote respectively
$ H^{i}(X^{\an},{\Bbb C}) $ and
$ {H}_{c}^{i}(X^{\an},{\Bbb C}) $.
A point of a variety means always a closed point, and
$ x \in  X $ means
$ x \in  X({\Bbb C}) $.

We denote the nearby and vanishing cycle functors [10] by
$ \psi $, $ \varphi  $,
and
$ \psi [-1] $,
$ \varphi [-1] $ by
$ {}^{p}\psi $,
$ {}^{p}\varphi  $,
because these preserve perverse sheaves.

\bigskip\bigskip
\centerline{{\bf 1.~Algebraic Gauss-Manin systems}}

\bigskip
\noindent
{\bf 1.1.}
Let
$ f : X \rightarrow  S $ be a morphism of smooth algebraic varieties.
Let
$ n = \dim X $.
We assume in this paper
$ \dim S = 1 $,
 $ f $ is not constant and
$ X, S $ are connected.
The Gauss-Manin systems
$ {\Cal G}_{f}^{i} $ are defined to be the cohomology sheaves
$ {\Cal H}^{i}{\Cal K}_{f} $ of the direct image
$ {\Cal K}_{f} := f_{+}({\Cal O}_{X}) $ of
$ {\Cal O}_{X} $ as algebraic left
$ {\Cal D} $-Modules.
We have
$$
\DR_{S}({\Cal G}_{f}^{i}) = {}^{p}R^{i+n}f_{*}{\Bbb C}_{X},
\tag 1.1.1
$$
where
$ {}^{p}R^{i}f_{*}{\Bbb C}_{X} = {}^{p}{\Cal H}^{i}(\bold{R}f_{*}{\Bbb C}_{X}) $ with
$ {}^{p}{\Cal H}^{i} : {D}_{c}^{b}(S,{\Bbb C}) \rightarrow   
\Perv(S,{\Bbb C}) $ the perverse cohomology functor in [3].

By assumption, we have a natural injective morphism
$ {\Cal O}_{S} \rightarrow  {\Cal G}_{f}^{1-n} $.
We define the reduced Gauss-Manin systems
$ {\tilde{\Cal G}}_{f}^{i} $ by
$$
{\tilde{\Cal G}}_{f}^{i} = 
\cases
\Coker({\Cal O}_{S} \rightarrow  {\Cal G}_{f}^{1-n})
&\text{for }i = 1 - n,\\ 
{\Cal G}_{f}^{i}
&\text{otherwise.}
\endcases
$$
Let
$ s \in  S $ with the natural inclusion
$ i_{s} : \{s\} \rightarrow  S $.
We choose and fix a local coordinate
$ t_{s} $ around
$ s $ such that
$ \{s\} = \{t_{s} = 0\} $.
For a complex of
$ {\Cal O}_{S} $-Modules (or
$ {\Cal O}_{S^{\an},s} $-modules)
$ M $,
we define
$$
{i}_{s}^{!}M =  \Cone(- t_{s} : M_{s} \rightarrow  M_{s})[-1].
\tag 1.1.2
$$

\proclaim{\bf 1.2.~Proposition}
We have short exact sequences
$$
0 \rightarrow  H^{1}{i}_{s}^{!}{\tilde{\Cal G}}_{f}^{i-1} \rightarrow  
{\tilde{H}}_{c}^{n-i}(X_{s},{\Bbb C})^{*} \rightarrow  
H^{0}{i}_{s}^{!}{\tilde{\Cal G}}_{f}^{i} \rightarrow  0,
\tag 1.2.1
$$
where
$ {\tilde{H}}_{c}^{n-i}(X_{s},{\Bbb C}) $ is as in the introduction, and
$ * $ denotes the dual vector space.
\endproclaim

\demo\nofrills {Proof.\usualspace}
Let
$ f_{+}({\Cal O}_{X})\tilde{\,\,} = C({\Cal O}_{S} \rightarrow  f_{+}
({\Cal O}_{X})) $.
Then 
$ {\Cal H}^{i}f_{+}({\Cal O}_{X})\tilde{\,\,} = {\tilde{\Cal G}}_{f}^{i} $, and
we have a spectral sequence
$$
{E}_{2}^{p,q} = H^{p}{i}_{s}^{!}{\tilde{\Cal G}}_{f}^{q} \Rightarrow  
H^{p+q}{i}_{s}^{!}f_{+}({\Cal O}_{X})\tilde{\,\,},
$$
degenerating at
$ E_{2} $ (because
$ {E}_{2}^{p,q} = 0 $ for
$ p \ne  0, 1) $.
So it is enough to show
$$
H^{i}{i}_{s}^{!}f_{+}({\Cal O}_{X})\tilde{\,\,} = {\tilde{H}}_{c}^{n-i}(X_{s},
{\Bbb C})^{*}.
$$
Then, using the distinguished triangle
$ \rightarrow  {\Cal O}_{S} \rightarrow  f_{+}({\Cal O}_{X}) \rightarrow  
f_{+}({\Cal O}_{X})\tilde{\,\,} \rightarrow  $ together with
$ {i}_{s}^{!}{\Cal O}_{S} = {\Bbb C}[-1],  $ we can reduce the assertion to
$$
H^{i}{i}_{s}^{!}f_{+}({\Cal O}_{X}) = {H}_{c}^{n-i}(X_{s},{\Bbb C})^{*}.
$$

Since the de Rham functor
$ \DR $ commutes with
$ {i}_{s}^{!} $ and the direct images, we have
$$
{i}_{s}^{!}f_{+}({\Cal O}_{X}) = \DR_{\{s\}}({i}_{s}^{!}f_{+}({\Cal O}_{X})) 
= {i}_{s}^{!}\bold{R}f_{*}\DR_{S}({\Cal O}_{X}) = 
{i}_{s}^{!}\bold{R}f_{*}({\Bbb C}_{X}[n]),
$$
because
$ \DR_{\{s\}} = id $ and
$ \DR_{S}({\Cal O}_{X}) = {\Bbb C}_{X}[n] $.
Since
$ {\Bbb C}_{X}[n] $ is self dual (i.e.,
$ {\Bbb D}({\Bbb C}_{X}[n]) = {\Bbb C}_{X}[n]),  $ the assertion follows 
from
$$
{\Bbb D}{i}_{s}^{!}\bold{R}f_{*}({\Bbb C}_{X}[n]) = {i}_{s}^{{}^*}f_{!}
{\Bbb D}({\Bbb C}_{X}[n]) = \bold{R}\Gamma _{c}(X_{s},{\Bbb C})[n],
$$
where
$ {\Bbb D} $ denotes the dual in the derived category of bounded 
complexes of
$ {\Bbb C} $-Modules with constructible cohomologies.
\enddemo

\proclaim{\bf 1.3.~Proposition}
Let
$ X = {\Bbb A}^{n} $, and
$ S = {\Bbb A}^{1} $.
Then
$ {\Cal G}_{f}^{i} $ has no nontrivial quotient with discrete support, and
$$
{}^{p}R^{i}f_{*}{\Bbb C}_{X} = (R^{i-1}f_{*}{\Bbb C}_{X})[1]
\quad\text{for any }i.
\tag 1.3.1
$$
If furthermore
$ f $ has at most isolated singularities including at infinity [31],
then
$$
\tilde{\Cal G}_{f}^{i} = 0 \quad\text{for }i \ne 0.
\tag 1.3.2
$$
\endproclaim

\demo\nofrills {Proof.\usualspace}
We have a long exact sequence
$$
\rightarrow H^{i+n-1}(X,{\Bbb C})
\rightarrow {G}_{f}^{i} \overset{\partial_{t}}\to\rightarrow
 {G}_{f}^{i} \rightarrow H^{i+n}(X,{\Bbb C}) \rightarrow,
$$
because
$ \bold{R}\Gamma(X,{\Bbb C})[n] = 
C(\partial_{t} : \bold{R}\Gamma(S,{\Cal K}_{f}) \rightarrow
\bold{R}\Gamma(S,{\Cal K}_{f})) $. 
See [4].
This implies the surjectivity of the action of
$ \partial_{t} $ on
$ {G}_{f}^{i} $, and the first assertion follows.

If
$ f $ has at most isolated singularities including at infinity
(or, more generally, if 
$ \supp {\Cal E}_{s} $ in (3.1) has discrete support), then we
see that
$ \varphi_{t_{s}} {}^{p}R^{i+n}f_{*}{\Bbb C}_{X} = 0 $ for
$ i \ne 0 $, using the commutativity of the vanishing cycle functors
with the direct image under a proper morphism
(because the vanishing cycles form a perverse sheaf with discrete support).
This implies that
$ {}^{p}R^{i+n}f_{*}{\Bbb C}_{X} $ is constant for
$ i \ne 0 $.
Then it vanishes for
$ i \ne 1-n, 0 $
by the above exact sequence, and 
$ {}^{p}R^{1}f_{*}{\Bbb C}_{X} = {\Bbb C}_{S}[1] $.
So the assertion follows from the Riemann-Hilbert correspondence.
\enddemo

\bigskip\bigskip
\centerline{\bf 2.~Algebraic Brieskorn modules}

\bigskip
\noindent
{\bf 2.1.}
Let
$ f : X \rightarrow  S $ and
$ {\Cal K}_{f} $ be as in (1.1).
Let
$ \omega _{S}\otimes _{{\Cal O}_{S}}{\Cal K}_{f} $ be the complex of 
right
$ {\Cal D}_{S} $-Modules corresponding to the direct image
$ {\Cal K}_{f} \,(= f_{+}({\Cal O}_{X})) $.
Then by definition of the direct image of right
$ {\Cal D} $-Modules, we have
$$
\omega _{S}\otimes _{{\Cal O}_{S}}{\Cal K}_{f} = 
\bold{R}f_{*}{\Cal C}_{f},\quad \omega _{S}
\otimes _{{\Cal O}_{S}}{\Cal G}_{f}^{i} = R^{i}f_{*}{\Cal C}_{f},
$$
where
$ {\Cal C}_{f} $ is the complex of {\it right}
$ f^{-1}{\Cal D}_{S} $-Modules such that
$$
{\Cal C}_{f}^{j} = {\Omega }_{X}^{j+n}\otimes _{f^{-1}{\Cal O}_{S}}
f^{-1}{\Cal D}_{S},
$$
and the differential is given by
$ ^{r}d_{f}(\omega \otimes P) = d\omega \otimes P + (f^{*}dt\wedge 
\omega )\otimes \partial _{t}P $ for
$ \omega  \in  {\Omega }_{X}^{j+n}, P \in  f^{-1}{\Cal D}_{S} $,
if
$ t $ is a local coordinate of
$ S $ and
$ \partial _{t} = \partial /\partial t $.

We define the filtration
$ F' $ on
$ {\Cal C}_{f} $ by
$ F'_{p}{\Cal C}_{f}^{j} = 0 $ for
$ p < -1 $ and
$$
\aligned
F'_{-1}{\Cal C}_{f}^{j} &= (f^{*}{\Omega }_{S}^{1}\wedge 
{\Omega }_{X}^{j+n-1})\otimes 1,
\\
F'_{p}{\Cal C}_{f}^{j} &= F'_{-1}{\Cal C}_{f}^{j} + 
{\Omega }_{X}^{j+n}\otimes f^{-1}F_{p+j}{\Cal D}_{S},
\endaligned
\tag 2.1.1
$$
where
$ F $ on
$ {\Cal D}_{S} $ is the filtration by the order of operators.
(This filtration
$ F' $ is different from the Hodge filtration
$ F $ in [27], and is useful only in the isolated singularity case.)
Let
$$
{\Cal C}_{f}^{(p)} = F'_{-p}{\Cal C}_{f},\quad {\Cal G}_{f}^{(p)} = 
{\omega }_{S}^{\vee }\otimes _{{\Cal O}_{S}}R^{0}f_{*}{\Cal C}_{f}^{(p)},
\tag 2.1.2
$$
where
$ {\omega }_{S}^{\vee } $ denotes the dual line bundle of
$ \omega _{S} $.

If
$ f $ is affine and
$ \dim  \Sing  f = 0 $,
then we have a natural morphism
$ {\Cal G}_{f,s}^{(0),\an} \rightarrow  {\Cal G}_{f,x}^{(0)} $ for
$ x\in \Sing X_{s},  $ where
$ {\Cal G}_{f,x}^{(0)} $ is the local analytic Brieskorn lattice of
$ f $ at
$ x $ [7].
Here
$ \omega _{S,s} $ is trivialized by using
$ t_{s} $ in (1.1).
In this case, we define
$$
{\Cal L}_{f,s}^{i} = 
\cases
\Ker({\Cal G}_{f,s}^{(0),\an} \rightarrow  \bigoplus _{x\in \Sing X_{s}} 
{\Cal G}_{f,x}^{(0)})
&\text{if  } i = 0,\\
{\tilde{\Cal G}}_{f,s}^{i,\an}
&\text{otherwise.}
\endcases
\tag 2.1.3
$$
See (1.1) for
$ {\tilde{\Cal G}}_{f}^{i} $.
We define also
$$
{\Cal G}_{f}^{\prime i} = 
\cases
{\Cal G}_{f}^{(0)}
&\text{if  } i = 0,\\
{\Cal G}_{f}^{i}
&\text{otherwise.}
\endcases
\tag 2.1.4
$$

\bigskip
\noindent
{\bf 2.2.}
{\it Remarks.} (i) If we choose a local coordinate
$ t $ (by shrinking
$ S) $,
then
$ \omega _{S} $ is trivialized by
$ dt $,
and it is well-known that right
$ {\Cal D}_{S} $-Modules are identified with left
$ {\Cal D}_{S} $-Modules by using the involution
$ * $ of
$ {\Cal D}_{S} $ defined by
$$
(PQ)^{*} = Q^{*}P^{*},\quad \partial _{t}{}^{*} = -\partial _{t},\,\,\,g^{*} 
= g\quad  \text{for }g \in  {\Cal O}_{S}.
$$
So we get an isomorphism
$ {\Cal K}_{f} = \bold{R}f_{*}{\Cal C}_{f} $ in the derived category of 
left
$ {\Cal D}_{S} $-Modules, and
$$
{\Cal C}_{f}^{j} = {\Omega }_{X}^{j+n}\otimes _{{\Bbb C}}
{\Bbb C}[\partial _{t}],
$$
with the differential
$ ^{l}d_{f}(\omega \otimes {\partial }_{t}^{j}) = d\omega \otimes 
{\partial }_{t}^{j} - (f^{*}dt\wedge \omega )\otimes {\partial 
}_{t}^{j+1} $ as is well-known.
The action of
$ f^{-1}{\Cal D}_{S} $ on
$ {\Cal C}_{f} $ is given by
$$
{\partial }_{t}^{i}(\omega \otimes {\partial }_{t}^{j}) = \omega \otimes 
{\partial }_{t}^{i+j},\quad g(\omega \otimes 1) = (f^{*}g)\omega 
\otimes 1\quad  \text{for }g \in  {\Cal O}_{S}.
$$

(ii) Let
$ X' = X \setminus   \Sing  f $.
Then, choosing a local coordinate
$ t $, we have an isomorphism of complexes
$$
f^{*}dt\wedge  : {\Omega }_{X'/S}^{\ssbull }[n-1] \rightarrow  F'_{p}
{\Cal C}_{f}|_{X'}\quad  \text{for }p \ge  -1.
\tag 2.2.1
$$

(iii) If
$ f $ is affine, we have
$$
R^{i}f_{*}F'_{p}{\Cal C}_{f} = {\Cal H}^{i}(f_{*}F'_{p}{\Cal C}_{f}).
\tag 2.2.2
$$
In particular,
$ R^{0}f_{*}F'_{0}{\Cal C}_{f} = f_{*}{\Omega }_{X}^{n}/(f^{*}dt\wedge 
df_{*}{\Omega }_{X}^{n-2}) $ locally (choosing a local coordinate
$ t) $,
because
$$
d(f_{*}(f^{*}dt\wedge {\Omega }_{X}^{n-2})) = d(f^{*}dt\wedge 
f_{*}{\Omega }_{X}^{n-2}) = f^{*}dt\wedge df_{*}{\Omega }_{X}^{n-2}
$$
by
$ f_{*}(\Im(f^{*}dt\wedge  : {\Omega }_{X}^{n-2} \rightarrow  {\Omega 
}_{X}^{n-1})) = \Im(f^{*}dt\wedge  : f_{*}{\Omega }_{X}^{n-2} 
\rightarrow  f_{*}{\Omega }_{X}^{n-1}) $.

\proclaim{\bf 2.3.~Lemma}
With the notation of 2.1, assume
$ X $ affine and
$ \dim  \Sing  f = 0 $.
Then the natural morphism
$$
R^{i}f_{*}F'_{p}{\Cal C}_{f} \rightarrow  \omega _{S}
\otimes _{{\Cal O}_{S}}{\Cal G}_{f}^{i}
\tag 2.3.1
$$
is an isomorphism for
$ p \ge  -1, i \ne  0 $,
and injective for
$ i = 0 $.
In particular, we have
$$
R^{i}f_{*}{\Cal C}_{f}^{(0)} = {\Cal G}_{f}^{\prime i}\quad 
 \text{for any }i.
\tag 2.3.2
$$
\endproclaim

\demo\nofrills {Proof.\usualspace}
Choosing a local coordinate
$ t $,
$ {\Gr}_{p}^{F'}{\Cal C}_{f} $ is locally isomorphic to
$ \tau '_{\ge -p}({\Omega }_{X}^{\ssbull }[n], f^{*}dt\wedge ) $ for
$ p \ge  0 $ 
where
$ \tau '_{\ge -p}K $ for a complex
$ (K, d) $ is defined by
$$
(\tau '_{\ge -p}K) ^{i} = \Coker \,d^{p-1}\,\,\, \text{if }
i = p,\,\,\,K^{i}\,\,\, \text{if }i > p,\,\,\, \text{and }0\,\,\, 
\text{otherwise } .
$$
Then
$ {\Gr}_{p}^{F'}{\Cal C}_{f} $ is quasi-isomorphic to
$ {\Omega }_{X/S}^{n} $ for
$ p \ge  0 $ by hypothesis, and the assertion follows from the long 
exact sequence
$$
\rightarrow  R^{i}f_{*}F'_{p-1}{\Cal C}_{f} \rightarrow  
R^{i}f_{*}F'_{p}{\Cal C}_{f} \rightarrow  
R^{i}f_{*}{\Gr}_{p}^{F'}{\Cal C}_{f} \rightarrow  
R^{i+1}f_{*}F'_{p-1}{\Cal C}_{f} \rightarrow  
$$
because
$ R^{i}f_{*}F'_{p-1}{\Cal C}_{f} = 0 $ for
$ i > 0 $ by the assumption on
$ f $.

\proclaim{\bf 2.4.~Theorem}
With the notation of (2.1), assume
$ f $ affine and
$ \dim  \Sing  f = 0 $.
Then we have short exact sequences
$$
0 \rightarrow  H^{1}{i}_{s}^{!}{\Cal L}_{f,s}^{i} \rightarrow  \tilde{H}^{i+n-1}
(X_{s}, {\Bbb C}) \rightarrow  H^{0}{i}_{s}^{!}{\Cal L}_{f,s}^{i+1} 
\rightarrow  0.
\tag 2.4.1
$$
\endproclaim

\demo\nofrills {Proof.\usualspace}
Let
$ \bar{f} : \bar{X} \rightarrow  S $ be a relative 
compactification of
$ f : X \rightarrow  S $ with
$ j : X \rightarrow  \bar{X} $ the inclusion morphism such that
$ \bar{X} $ is smooth and
$ D := \bar{X} \setminus  X $ is a divisor.
We define
$ {\Cal C}_{\bar{f}}^{(0)} $ on
$ \bar{X} $ as in (2.1).
Let
$ {\Cal C}_{\bar{f}}^{(0)}(*D) $ be the localization of
$ {\Cal C}_{\bar{f}}^{(0)} $ by a local defining equation of
$ D $.
Then
$ {\Cal C}_{\bar{f}}^{(0)}(*D) = j_{*}{\Cal C}_{f}^{(0)} $,
and we have a natural surjective morphism
$$
{\Cal C}_{\bar{f}}^{(0)}(*D)^{\an} \rightarrow  
\bigoplus _{x\in \Sing f} {\Cal G}_{f,x}^{(0)},
$$
where
$ {\Cal G}_{f,x}^{(0)} $ is viewed as a sheaf on
$ \bar{X} $ with support
$ \{x\} $.
Indeed, by the theory of Gauss-Manin connection [7], we have
$$
{\Cal H}^{i}({\Cal C}_{f}^{(0)}|{X}_{s}^{\an}) = 
\cases
{\Bbb C}\{t_{s}\}\bigotimes _{{\Bbb C}}{\Bbb C}_{X_{s}} 
& \text{if }i = 1 - n,\\
\bigoplus _{x \in  \Sing X_{s}} {\Cal G}_{f,x}^{(0)}
&\text{if }i = 0,\\
0
&\text{otherwise.}
\endcases
$$

Let
$ {\Cal C}'_{f} =  \Cone({\Cal C}_{\bar{f}}^{(0)}(*D)^{\an} 
\rightarrow  \bigoplus _{x\in \Sing f} {\Cal G}_{f,x}^{(0)})[-1] $.
Then
$$
{\Cal C}'_{f}|{X}_{s}^{\an} ={\Bbb C}\{t_{s}\}\otimes _{{\Bbb C}}
{\Bbb C}_{X_{s}}[n-1].
\tag 2.4.2
$$
Since
$ f $ is affine and
$ t_{s}{\Cal G}_{f,x}^{(0)} \supset  {\partial }_{{t}_{s}}^{-k}
{\Cal G}_{f,x}^{(0)} $ for
$ k \gg 0 $,
Nakayama's lemma implies
$$
{\Cal G}_{f,s}^{(0),\an} \rightarrow  \bigoplus _{x \in  \Sing X_{s}} 
{\Cal G}_{f,x}^{(0)}\,\,\, \text{is surjective } .
\tag 2.4.3
$$
So we get
$$
{\Cal L}_{f,s}^{i} = (R^{i}\bar{f}_{*}{\Cal C}'_{f})_{s}\quad  \text{for }i 
\ne  1 - n.
$$
For
$ i = 1 - n $,
we have a natural injective morphism
$ {\Cal O}_{S^{\an}} \rightarrow  R^{1-n}\bar{f}_{*}{\Cal C}'_{f} $,
and
$$
{\Cal L}_{f,s}^{1-n} = \Coker({\Cal O}_{S^{\an}} \rightarrow  R^{1-n}
\bar{f}_{*}{\Cal C}'_{f})_{s}.
$$

Let
$ {i}_{s}^{{}^*}{\Cal C}'_{f} = \Cone(t_{s} : {\Cal C}'_{f} \rightarrow  
{\Cal C}'_{f})|{\bar{X}}_{s}^{\an} $.
Let
$ \bar{j}_{s} : X_{s} \rightarrow  \bar{X}_{s} $ denote the 
inclusion morphism.
Then it is enough to show
$$
{i}_{s}^{{}^*}{\Cal C}'_{f} = \bold{R}(\bar{j}_{s})_{*}{\Bbb C}_{X}
[n-1],
\tag 2.4.4
$$
using a spectral sequence similar to that in the proof of (1.2).
But
$$
{i}_{s}^{{}^*}{\Cal C}'_{f}|{X}_{s}^{\an} = {\Bbb C}_{X}[n-1]
$$
by (2.4.2), and the assertion is reduced to showing a natural 
quasi-isomorphism
$$
{i}_{s}^{{}^*}{\Cal C}'_{f} \simto 
\bold{R}(\bar{j}_{s})_{*}(\bar{j}_{s})^{*}{i}_{s}^{{}^*}
{\Cal C}'_{f}.
\tag 2.4.5
$$

Since the assertion is restricted to a neighborhood of
$ D $,
we may restrict to
$ \bar{X}' := \bar{X} \setminus   \Sing  f $,
and replace
$ {\Cal C}'_{f} $ with
$ {\Cal C}_{\bar{f}}^{(0)}(*D)^{\an} $.
Let
$ X' = X \setminus   \Sing  f, X'_{s} = X' \cap  X_{s} $ with the inclusion 
morphisms
$ j' : X' \rightarrow  \bar{X}', j'_{s} : X'_{s} \rightarrow  
\bar{X}'_{s} $.
Then
$ {\Cal C}_{\bar{f}}^{(0)}(*D)|\bar{X}' = 
j'_{*}{\Omega }_{X'/S}^{\ssbull } $ by (2.2.1).
So (2.4.5) is verified by applying the functor
$ \an $ to the distinguished triangle
$$
\rightarrow  j'_{*}{\Omega }_{X'/S}^{\ssbull } 
\overset{{t}_{s}}\to\rightarrow 
j'_{*}{\Omega }_{X'/S}^{\ssbull } \rightarrow  
(j'_{s})_{*}{\Omega }_{X'_{s}}^{\ssbull } \rightarrow .
$$
This completes the proof of (2.4).
\enddemo

\proclaim{\bf 2.5.~Corollary}
With the notation of (2.1) and the assumption of (2.4), let
$$
{N}_{s}^{\prime i} = \dim \Ker(t_{s} : {\Cal G}_{f}^{\prime i} 
\rightarrow  {\Cal G}_{f}^{\prime i}),
\quad {R}_{s}^{\prime i} = \dim \Coker(t_{s} : {\Cal G}_{f}^{\prime i} 
\rightarrow  {\Cal G}_{f}^{\prime i}),
$$
for
$ s \in  S $,
and let
$ \mu _{x} $ denote the Milnor number of
$ f $ at
$ x \in  X $.
Then
$$
{N}_{s}^{\prime i+1} + {R}_{s}^{\prime i} = \dim H^{i+n-1}(X_{s}, 
{\Bbb C}) + \delta _{i,0} \sum _{x\in \Sing X_{s}} \mu _{x}.
$$
\endproclaim

\demo\nofrills {Proof.\usualspace}
This follows from (2.4) because the local analytic Brieskorn lattice
$ {\Cal G}_{f,x}^{(0)} $ is a free
$ {\Bbb C}\{t\} $-module of rank
$ \mu _{x} $ by [30].
\enddemo

\bigskip
\noindent
{\bf 2.6.}
{\it Remark.}
Let
$ s \in  \Delta  $ and
$ s' \in  U $ such that
$ s' $ is sufficiently near
$ s $.
Then there does not exist a natural morphism
$$
\iota  : H^{n-1}(X_{s},{\Bbb C}) \rightarrow  H^{n-1}(X_{s'},{\Bbb C})
$$
making the following diagram commutative :
$$
\CD
(R^{n-1}f_{*}{\Bbb C}_{X})_{s}
@<{\sim}<<
H^{n-1}(f^{-1}(D_{s}),{\Bbb C})
\\
@VVV @VVV
\\
H^{n-1}(X_{s},{\Bbb C})   @>\iota>>    H^{n-1}(X_{s'},{\Bbb C})
\endCD
$$
where the vertical morphisms are natural morphisms, and
$ D_{s} $ is a sufficiently small open disk with center
$ s  $ such that
$ s' \in  D_{s} \setminus \{s\} $.
Indeed, the right vertical morphism of the diagram is injective because
$ (R^{n-1}f_{*}{\Bbb C}_{X})[1] $ is a perverse sheaf.
But the left vertical morphism is not bijective if
$ n = 2 $ and
$ \nu _{s} \ne  0 $ due to Th.~3 of [1].
There exist examples such that
$ \dim (R^{n-1}f_{*}{\Bbb C}_{X})_{s} = \dim H^{n-1}(X_{s},{\Bbb C}) $ with
$ n = 2 $ and
$ \nu _{s} \ne  0 $ (e.g.
$ f = x^{4}y^{2} + 2x^{2}y + xy^{2} $ and
$ s = -1 $).

\bigskip\bigskip
\centerline{{\bf 3.~Vanishing cycles}}

\bigskip
\noindent
{\bf 3.1.}
With the notation of (1.1), let
$ \bar{f} : \bar{X} \rightarrow  S $ be a relative 
compactification of
$ f $ with 
$ j : X \rightarrow  \bar{X} $ the open immersion such that
$ \bar{f}j = f $.
(Here
$ \bar{f} $ is proper, but
$ \bar{X} $ may be singular.)
For
$ s \in  \Delta  := S \setminus  U $,
let
$$
{\Cal E}_{s} = \varphi _{\bar{f}^{*}t_{s}} \bold{R}j_{*}
{\Bbb C}_{X}[n-1],\quad 
{\Cal E}_{s,\infty } = {\Cal E}_{s}|\bar{X}_{s}\setminus X_{s},
\quad
{\Cal E}_{s,\fin} = {\Cal E}_{s}|X_{s}, 
$$
where
$ \varphi  $ denotes the vanishing cycle functor [10], and
$ \bar{X}_{s} = \bar{f}^{-1}(s),  $ etc.
We define
$$
\aligned
E_{s} = \bold{H}^{0}({\bar{X}}_{s}^{\an}, {\Cal E}_{s}),
&\quad
E'_{s} = \bold{H}^{0}({\bar{X}}_{s}^{\an}, {\Cal E}_{s,\infty }),
\\
E''_{s} = \bold{H}^{0}({X}_{s}^{\an}, {\Cal E}_{s,\fin}),
&\quad
E''_{x} = ({\Cal H}^{0}{\Cal E}_{s})_{x},
\endaligned
$$
for
$ x\in \Sing X_{s} $.
These vector spaces have naturally the monodromy
$ T $ associated with the functor
$ \varphi  $.
Since
$  \supp \, {\Cal E}_{s,\fin} =  \Sing  X_{s} $,
we have
$$
E_{s} = E'_{s} \oplus  E''_{s},\quad E''_{s}= \bigoplus _{f(x) =s} E''_{x}.
$$

Let
$ E_{s}^{\lambda} = \Ker(T_{ss}-\lambda|E_{s}) $ (similarly for
$ E_{s}^{\prime \lambda} $,
$ E_{s}^{\prime \prime \lambda} $), where
$ T_{ss} $ is the semisimple part of the monodromy
$ T $.
Let
$$
\nu_{s}^{\lambda} = \dim E_{s}^{\prime \lambda},
\quad
\mu_{s}^{\lambda} = \dim E_{s}^{\prime\prime \lambda},
\tag 3.1.1
$$
and
$ \nu_{s} = \sum_{\lambda} \nu_{s}^{\lambda} $,
$ \nu = \sum_{s} \nu_{s} $ (similarly for
$ \mu $).

We define
$ P_{s}(E,T) $ for
$ s \in  S $ and for a finite dimensional  
$ {\Bbb C} $-vector space
$ E $ with a quasi-unipotent automorphism
$ T $ as follows.
Let
$ t_{s} $ be as in (1.1), and we will write
$ \partial _{t} $ for
$ \partial _{t_{s}} $ to simplify the notation.
Let
$ M_{s}(E,T) = E\otimes _{{\Bbb C}}{\Bbb C}[\partial _{t},{\partial 
}_{t}^{-1}] $ with action of
$ {\partial }_{t}^{i} \,(i \in  {\Bbb Z}) $ and
$ t_{s} $ defined by
$$
\aligned
{\partial }_{t}^{i}(e\otimes {\partial }_{t}^{j}) 
&= e\otimes {\partial }_{t}^{i+j},
\\
t_{s}(e\otimes {\partial }_{t}^{j}) &= (1 - j)e\otimes 
{\partial }_{t}^{j-1} - (2\pi i)^{-1}(\log T)e\otimes {\partial }_{t}^{j-1}.
\endaligned
$$
Then
$ {\partial }_{t}^{i}t - t_{s}{\partial }_{t}^{i} = i{\partial }_{t}^{i-1} $ 
and
$ t_{s}\partial _{t} + j $ on
$ E\otimes {\partial }_{t}^{j} $ is
$ - (2\pi i)^{-1}\log T\otimes id $.
Here
$ \log  T $ is chosen so that the eigenvalues of
$ (2\pi i)^{-1}\log T $ are contained in
$ [0,1) $.
We define
$$
\aligned
M_{s}(E,T)^{>0} 
&= E\otimes _{{\Bbb C}}{\Bbb C}[{\partial }_{t}^{-1}],
\\
P_{s}(E,T) &= M_{s}(E,T)/M_{s}(E,T)^{>0}
\endaligned
$$
They are
$ {\Bbb C}[t_{s}] $-modules with action of
$ {\partial }_{t}^{-1} $.

By a canonical splitting of the filtration
$ V $ in the one variable case (see [26, 1.5]) together with the 
commutativity of the de Rham functor with the vanishing cycle functor
[18], [20] (see also [25, 3.4.12]), we get canonical isomorphisms
$$
{G}_{f}/{G}_{f}^{>0} = \bigoplus _{s\in \Delta } P_{s}(E_{s},T),
\quad
{\Cal G}_{f,x}/{\Cal G}_{f,x}^{>0} = P_{f(x)}(E''_{x},T),
\tag 3.1.2
$$
because
$$
E_{s} = {}^{p}\varphi _{t_{s}}({}^{p}R^{n}f_{*}{\Bbb C}_{X}) = 
{}^{p}\varphi _{t_{s}}\DR_{S}({\Cal G}_{f}^{\an}),
\quad 
E''_{x} = {}^{p}\varphi _{t_{s}}\DR_{S'}({\Cal G}_{f,x}),
$$
where
$ S' $ is an analytic open neighborhood of
$ s $,
and
$ {\Cal G}_{f,x} $ denotes also its coherent extension to
$ S' $.

\proclaim{\bf 3.2.~Theorem}
Let  $ f : X \rightarrow  S $ be as in (1.1) and assume
$ X, S $ affine and
$ \dim  \Sing  f = 0 $.
Then, with the notation of (3.1), we have the commutative diagram of 
exact sequences in (0.5), where
$ {\Bbb C}[t], c, t - c $ are replaced respectively by
$ \Gamma (S, {\Cal O}_{S}), s $ and
$ t_{s} $.
Furthermore, this holds with
$ {G}_{f}^{(0)} $,
$ {G}_{f}^{(0),>0} $,
$ {\Cal G}_{f,x}^{(0)} $ replaced by
$ {G}_{f}^{(-1)} $,
$ {G}_{f}^{(-1),>0} $,
$ {\Cal G}_{f,x}^{(-1)} $.
If there exists a nowhere vanishing vector field
$ v $ on
$ S $,
and the action of
$ v $ is bijective on
$ G_{f} $,
then the morphisms of the diagram are compatible with the action of
$ v^{-1} $.
\endproclaim

\demo\nofrills {Proof.\usualspace}
We prove the assertion for
$ {G}_{f}^{(0)} $.
The argument is similar for
$ {G}_{f}^{(-1)} $.

The exactness of the middle and right columns follows from (3.1.2).
The morphism between the two columns is defined by 
taking the restriction to the local Milnor fibration at
$ x $.
So the morphisms of the diagram are naturally defined except for the 
surjective morphism of the left column,
 but it is induced by the other morphisms using the commutativity and 
the exactness.

By  $ \GAGA,  $ we may replace
$ G_{f} $,
$ {G}_{f}^{(0)} $,
$ {G}_{f}^{>0} $ by
$ {\Cal G}_{f}^{\an} $,
$ {\Cal G}_{f}^{(0),\an} $,
$ {\Cal G}_{f}^{>0,\an} $ respectively, and then restrict to the stalk at
$ s \in  S $,
where
$ P_{s}(E_{s},T),  $ etc.~are viewed as sheaves with support
$ \{s\} $.
So it is enough to show the exactness of the middle row 
and the surjectivity of the second morphism in the upper row, 
because the bottom row splits, and is exact.

We consider a morphism
$$
{\Cal G}_{f,s}^{\an} \rightarrow  \bigoplus _{x \in  \Sing X_{s}} 
{\Cal G}_{f,x},
\tag 3.2.1
$$
which is surjective by (2.4.3).
The source has the filtration
$ F' $ defined in (2.1) such that
$$
F'_{0}{\Cal G}_{f,s}^{\an} = {\Cal G}_{f,s}^{(0),\an},\quad 
{\Gr}_{p}^{F'}{\Cal G}_{f,s}^{\an} = (f_{*}{\Omega }_{X/S}^{n})_{s}\quad  
\text{for }p > 0.
\tag 3.2.2
$$
See (2.2.2).
On
$ {\Cal G}_{f,x} \,(= {\Cal H}^{0}{\Cal C}_{f,x}^{\an}) $,
 $ F' $ induces the filtration
$ F' $ satisfying a similar property where
$ {\Cal G}_{f,s}^{\an}$,
$ {\Cal G}_{f,s}^{(0),\an} $,
$ (f_{*}{\Omega }_{X/S}^{n})_{s} $ are replaced respectively by
$ {\Cal G}_{f,x} $,
$ {\Cal G}_{f,x}^{(0)} $,
$ {\Omega }_{X/S,x}^{n} $.
So (3.2.1) induces an isomorphism
$$
{\Gr}_{p}^{F'}{\Cal G}_{f,s}^{\an} \rightarrow  \bigoplus _{x \in  \Sing 
X_{s}} {\Gr}_{p}^{F'}{\Cal G}_{f,x}\quad  \text{for }p > 0,
$$
and hence
$ {\Cal G}_{f,s}^{\an}/{\Cal G}_{f,s}^{(0),\an} \rightarrow  
\bigoplus _{x \in  \Sing X_{s}} {\Cal G}_{f,x}/{\Cal G}_{f,x}^{(0)} $ is an 
isomorphism.
This shows the exactness of the middle row.

Now it remains to show the surjectivity of
$$
{\Cal G}_{f,s}^{>0,\an} \rightarrow  \bigoplus _{x \in  \Sing X_{s}} 
{\Cal G}_{f,x}^{>0}.
$$
But (3.2.1) is surjective, and the assertion follows from the property of 
the Deligne extension (or the filtration
$ V) $.
See for example [24, 3.1.5].

The last assertion on the action on
$ v^{-1} $ is clear.
\enddemo

\noindent
{\bf 3.3.}
{\it Remark.} Assume
$ f $ is affine, and has at most isolated singularities including at infinity.
Let
$ {G}_{f}^{(0)} $ be as in (2.1).
As a corollary of the results of Sabbah [33],
$ f $ is cohomologically tame in his sense if and only if
$ {G}_{f}^{(0)} $ is coherent over
$ {\Cal O}_{S} $.
For the ``if'' part it is sufficient to show
$ {G}_{f}^{(0)} \subset  {G}_{f}^{>0} $,
and this is easily verified by using the filtration of Kashiwara [18] and 
Malgrange [20].
See loc.~cit.
For the converse, it is sufficient, however, to note that, for a good 
filtration
$ F $ of a (regular) holonomic
$ {\Cal D}_{S} $-Module
$ M $,
we have for each
$ s \in  S $
$$
\dim {\Gr}_{p}^{F}M_{s} = \dim {}^{p}\varphi _{t_{s}} \DR_{S}(M)\quad  
\text{for }p \gg  0.
\tag 3.3.1
$$
Indeed, if
$ (M,F) = ({\Cal G}_{f},F') $,
we have  $ \dim {\Gr}_{p}^{F}M_{s} = \mu _{s} $ for
$ p \gg  0 $,
and the above equality gives
$ \nu _{s} = 0 $.
For the proof of (3.3.1), note that
$ \dim {\Gr}_{p}^{F}M_{s} $ for
$ p \gg  0 $ coincides with the multiplicity, and is independent of the 
choice of the good filtration
$ F $.
(In the regular holonomic case, we can take
$ F_{p} = V^{\alpha -p} \, (p > 0) $ for some
$ \alpha  \in  {\Bbb Q} $,
where
$ V $ denotes the filtration of Kashiwara [18] and Malgrange [20]
in the one variable case.)

\bigskip
\noindent
{\bf 3.4.} {\it Proof of Corollary }(0.7).
Let
$ S $ be a smooth analytic curve, and
$ \Delta  $ a discrete subset.
Put
$ U  = S \setminus  \Delta  $.
Let
$ L $ be a
$ {\Bbb C} $-local system on
$ U $ with quasi-unipotent local monodromies around any
$ s \in  \Delta  $,
and
$ {\Cal L} $ the meromorphic Deligne extension.
That is,
$ {\Cal L} $ is a regular holonomic
$ {\Cal D}_{S} $-Module such that
$ {\Cal L}|_{U} = L\otimes _{{\Bbb C}}{\Cal O}_{U} $ and for any
$ s \in  \Delta  $,
the action of local equation of
$ \{s\} $ is bijective on
$ {\Cal L}_{s} $.
Let
$ {\Cal L}^{(0)} $ be a coherent
$ {\Cal O}_{S} $-submodule of
$ {\Cal L} $ such that
$ {\Cal L}^{(0)}|_{U} = {\Cal L}|_{U} $.
We define
$ \ord_{s} {\Cal L}^{(0)} $, the order of
$ {\Cal L}^{(0)} $ at
$ s \in  \Delta  $, as follows.

Let
$ S' $ be a sufficiently small open disk around
$ s $ with a coordinate
$ t $,
and
$ \{\omega _{i}\} $ an
$ {\Cal O}_{S'} $-basis of
$ {\Cal L}|_{S'} $.
Let
$ L^{*} $ be the dual local system of
$ L $,
and
$ \{e_{j}\} $ a basis of horizontal (multivalued) sections of
$ L^{*}|_{S'\setminus \{s\}} $.
Then
$ \langle e_{j},\omega _{i}\rangle $ is a Nilson class function on
$ S' \setminus  \{s\} $,
and the determinant has the asymptotic expansion
$$
\det\langle e_{j},\omega _{i}\rangle =
 C t^{\alpha } +  \text{higher terms } 
$$
for
$ C \in  {\Bbb C}^{*} $ and
$ \alpha  \in  {\Bbb Q} $, because 
$ \det L $ has local monodromies with finite order.
This
$ \alpha  $ is independent of the choice of the bases, and we define
$ \ord_{s} {\Cal L}^{(0)} = \alpha  $.
Then we get (0.7) from the following observations.

\bigskip

(i) If
$ L $ has a finite filtration
$ G $,
then it induces the filtration
$ G $ on
$ {\Cal L} $ and
$ {\Cal L}^{(0)} $,
and
$$
\ord_{s} {\Cal L}^{(0)} = 
\prod _{j} \ord_{s} {\Gr}_{j}^{G}{\Cal L}^{(0)}.
$$

(ii) If
$ {\Cal L}^{(0)} $ is the Deligne extension such that the eigenvalues of 
the residues are contained in
$ (\alpha , \alpha +1] $ (resp.
$ [\alpha , \alpha +1)) $,
then
$$
\ord_{s} {\Cal L}^{(0)} = \sum\nolimits\limits_{j=1}^{r} {\alpha }_{j},
$$
where
$ r =  \rank \, L $,
and the
$ \alpha _{j} $ are rational numbers contained in
$ (\alpha , \alpha +1] $ (resp.
$ [\alpha , \alpha  +1)) $ such that
$ \exp(-2\pi i\alpha _{j}) $ are the eigenvalues of the local monodromy 
of
$ L $ around
$ s $ (with multiplicity).

(iii) If
$ {\Cal L}^{(0)} $ is a coherent extension of the local analytic Brieskorn 
lattice
$ {\Cal G}_{f,x}^{(0)} $ (resp.
$ {\Cal G}_{f,x}^{(-1)} $),
then
$$
\ord_{s} {\Cal L}^{(0)} = (n - 2)\mu _{x}/2
\quad\text{(resp. }n\mu_{x}/2).
$$
This is well-known after [36].

(iv) If
$ {\Cal L}^{(0)} $ is
$ {\Cal G}_{f,s}^{(0),>0,\an} $ 
with the assumption of (3.2), then
$$
\ord_{s} {\Cal L}^{(0)} = 
(- {\nu }_{s}^{\ne 1} + \sum _{x\in \Sing X_{s}} (n - 2)\mu _{x})/2.
$$
where
$ {\nu }_{s}^{\ne 1} = \sum _{\lambda \ne 1} {\nu }_{s}^{\lambda } $.
This follows from the above observations together with the short exact 
sequence
$$
0 \rightarrow   \Kernel \rightarrow  {\Cal G}_{f,s}^{(0),>0,\an} 
\rightarrow  \bigoplus _{x \in  \Sing X_{s}} {\Cal G}_{f,x}^{(0)} 
\rightarrow  0,
$$
where
$  \Kernel = \Ker({\Cal G}_{f,s}^{>0,\an} \rightarrow  
\bigoplus _{x \in  \Sing X_{s}} {\Cal G}_{f,x}^{>0}) $ by (3.2).
(The argument is similar for
$ {\Cal G}_{f,s}^{(-1),>0,\an} $ with
$ (n-2)\mu_{x} $ replaced by
$ n\mu_{x} $.)

(v) If
$ L $ and
$ e_{j} $ are defined over
$ {\Bbb Z} $,
then
$ \det\langle e_{j},\omega _{i}\rangle^{2} $ is a meromorphic function on
$ S $.
If furthermore
$ S $,
$ {\Cal L} $,
$ {\Cal L}^{(0)} $ and
$ \omega_{i} $ are algebraic, and
$ {\Cal L} $ is also regular at 
infinity, then
$ \det\langle e_{j},\omega _{i}\rangle^{2} $ is a rational function on
$ S $.

\proclaim{\bf 3.5.~Proposition}
With the notation and assumption of (3.2), let
$ F'_{p}{G}_{f} = \partial_{t}^{p}{G}_{f}^{(0)} $,
$ F_{p}{\Cal G}_{f,x} = \partial_{t}^{p}{\Cal G}_{f,x}^{(0)} $, and let
$ V_{s} $ denote the filtration by the eigenvalue of the action of 
$ \partial_{t_{s}}t_{s} $ indexed by
$ {\Bbb Q} $ (i.e.
$ \partial_{t_{s}}t_{s} - \alpha $ is nilpotent on
$ {\Gr}_{V_{s}}^{\alpha} $).
See [9], [18], [20], etc.
Then
$$
{\Gr}_{p}^{F'}{\Gr}_{V_{s}}^{\alpha}{G}_{f} = 
\bigoplus_{x\in \Sing X_{s}}{\Gr}_{p}^{F}{\Gr}_{V_{s}}^{\alpha}{\Cal G}_{f,x}
$$
for 
$ p \ge 0 $ and 
$ 0 < \alpha \le 1 $.
In particular,
$ {G}_{f}^{(0)} $ together with the filtration
$ V_{s} $ gives the sums of the Hodge numbers of the local Milnor fibers.
\endproclaim

\demo\nofrills {Proof.\usualspace}
Let
$ F' $ denote also the corresponding filtration on
$ {\Cal G}_{f} $,
$ {\Cal G}_{f}^{\an} $.
Then we have a canonical isomorphism
$$
{\Gr}_{p}^{F'}{\Gr}_{V_{s}}^{\alpha}{G}_{f} 
= {\Gr}_{p}^{F'}{\Gr}_{V_{s}}^{\alpha}{\Cal G}_{f,s}^{\an}.
$$
by the exactness of the functors involved.
We have furthermore
$$
{\Gr}_{p}^{F'}{\Cal G}_{f,s}^{\an} 
= \bigoplus_{x\in \Sing X_{s}}{\Gr}_{p}^{F}{\Cal G}_{f,x},
$$
because 
$ {\Gr}_{p}^{F'}{G}_{f}
= \bigoplus_{s}{\Gr}_{p}^{F'}{\Cal G}_{f,s}^{\an} $.
So it is enough to show the bistrict surjectivity of
$$
r_{s} : ({\Cal G}_{f,s}^{\an};F',V_{s}) \rightarrow
\bigoplus_{x\in \Sing X_{s}}({\Cal G}_{f,x};F,V_{s}),
$$
where the filtration
$ F' $ is restricted to the index
$ p \ge -1 $.
Since the strict surjectivity is separately clear,
the assertion is reduced to the compatibility of three submodules
$ \Ker \, r_{s} $,
$ F'_{p} $,
$ V_{s}^{\alpha} $ of
$ {\Cal G}_{f,s}^{\an} $ due to [24, 1.2.14].
But this is clear because
$ \Ker\, r_{s} \subset F'_{p} $ for
$ p \ge -1 $ by (3.2).
The last assertion follows from [35] (see also [22], [27], [29], etc.)
\enddemo

\noindent
{\it Remark.}
We can also show for
$ \alpha \ne 0 $
$$
\dim {\Gr}_{V_{s}}^{\alpha+1}{G}_{f}^{(0)} -
\dim {\Gr}_{V_{s}}^{\alpha}{G}_{f}^{(0)} =
\sum_{x\in \Sing X_{s}}
(\dim {\Gr}_{V_{s}}^{\alpha+1}{\Cal G}_{f,x}^{(0)} -
\dim {\Gr}_{V_{s}}^{\alpha}{\Cal G}_{f,x}^{(0)} ),
$$
because
$ \dim {\Gr}_{V_{s}}^{\alpha+1}{\Cal L}_{f,s} =
\dim {\Gr}_{V_{s}}^{\alpha}{\Cal L}_{f,s} $ for
$ \alpha \ne 0 $ with the notation of (0.4).
This is closer to Varchenko's formulation [35].
However, it gives only the Hodge numbers
$ \dim {\Gr}_{F}^{p}H^{n-1}(X_{x},{\Bbb C})_{\lambda} $ for
$ p \ne n-1 $ or
$ \lambda \ne 1 $ (because
$ \alpha \ne 0 $), where 
$ X_{x} $ denotes the local Milnor fiber at
$ x \in \Sing X_{s} $, and the index
$ \lambda $ means the dimension of the
$ \lambda $-eigenspace by the monodromy.
To get the Hodge number for 
$ p = n - 1 $ and
$ \lambda = 1 $, we have to use also
$ G_{f}^{(-1)} $.

\bigskip\noindent
{\bf 3.6.}
{\it Relation between the numerical invariants.}
Let
$ X = {\Bbb A}^{n} $, and
$ S = {\Bbb A}^{1} $.
Assume
$ f $ has at most isolated singularities including at infinity [31] (or, 
more generally, the
$ \supp {\Cal E}_{s} $ are discrete.)
Let  $ \mu _{s}, \nu _{s}, \mu , \nu  $ be as in (3.1), and
$ m =  \rank \, L|_{U} \,(= \dim H^{n-1}(X_{s}, {\Bbb C}) $ for a general
$ s) $.
Then we can easily show that the vanishing of
$ H^{i}(S^{\an},L) $ for any
$ i $ implies
$$
m = \mu  + \nu ,
\tag 3.6.1
$$
which was first obtained by [31], [33].
Here apparently different definitions of
$ \nu $ are used, but they are all equivalent, because we can also show
$$
\chi (X_{s}) - \chi (X_{s'}) = (-1)^{n}(\nu _{s} + \mu _{s})
\tag 3.6.2
$$
for
$ s, s' \in  S $ such that
$ s' \notin  \Delta  $, and the corresponding formula is proved in loc. cit.

For
$ s \in  S $,
let
$ {\rho }_{s}^{\lambda }, {\rho }_{s}^{\prime \lambda } $ be the number 
of Jordan blocks of
$ T $ (i.e. the minimal number of generators over
$ {\Bbb C}[T]) $ on
$ {E}_{s}^{\lambda }, {E}_{s}^{\prime \lambda } $ respectively.
We define
$$
\aligned
R'_{s} = \dim \Coker(t_{s} : {G}_{f}^{(0)} \rightarrow  {G}_{f}^{(0)}),
\quad 
&R_{s} = \dim \Coker(t_{s} : G_{f} \rightarrow  G_{f}),
\\
N'_{s} = \dim \Ker(t_{s} : {G}_{f}^{(0)} \rightarrow  {G}_{f}^{(0)}),
\quad 
&N_{s} = \dim \Ker(t_{s} : G_{f} \rightarrow  G_{f}).
\endaligned
$$
Let
$ T_{s} $ denote the local monodromy of the local system
$ L|_{U} (= R^{n-1}f_{*}{\Bbb C}_{X}|_{U}) $ around
$ s \in  \Delta  $.
Let
$ \tilde{h}^{i}(X_{s}) = \dim \tilde{H}^{i}(X_{s},{\Bbb C}) $,
$ {\tilde{h}}_{c}^{i}(X_{s}) = \dim {\tilde{H}}_{c}^{i}(X_{s},{\Bbb C}) $.
Then the relation between
$ \tilde{h}^{i}(X_{s}) $,
$ {\tilde{h}}_{c}^{i}(X_{s}) $,
$ R'_{s} $,
$ R_{s} $,
$ N'_{s} $,
$ N_{s} $,
$ \nu _{s} $,
$ \mu _{s} $,
$ \mu _{x} $,
$ {\rho }_{s}^{\prime 1} $,
$ {\rho }_{s}^{1} $,
$ T_{s} $ and
$ m $ is summarized as follows.
$$
\aligned
\tilde{h}^{i}(X_{s}) = 0\,\,\, \text{if }i \ne  n - 2, n - 1,\quad 
&{\tilde{h}}_{c}^{i}(X_{s}) = 0\,\,\, \text{if }i \ne  n - 1, n,
\\
R'_{s} = \tilde{h}^{n-1}(X_{s}) + \sum _{x\in \Sing X_{s}} \mu _{x},\quad 
&R_{s} = {\tilde{h}}_{c}^{n-1}(X_{s}) = \dim \Coker(T_{s} - id),
\\
N'_{s} = \tilde{h}^{n-2}(X_{s}) \le  {\rho }_{s}^{\prime 1},\quad 
&N_{s} = {\tilde{h}}_{c}^{n}(X_{s}) \le  {\rho }_{s}^{1},
\\
R'_{s} - N'_{s} = m - \nu _{s},\quad 
&R_{s} - N_{s} = m - \mu _{s} - \nu _{s}.
\endaligned
$$
In particular, we can calculate the dimension of the cohomology 
(with compact support) of any fiber using the action of
$ t $ on
$ {G}_{f}^{(0)}, G_{f} $ together with the Milnor numbers
$ \mu _{x} $.
Most of the relations follow from (1.2--3), (2.4) and (3.6.1).
The formula involving
$ \Coker(T_{s}-id) $ is closely related to Th. 1 of [1],
and the vanishing of the Betti numbers to [34].
See also [12], [13], [14], [31], [33], etc.

In the case
$ n = 2 $, we can show that the monodromy
$ T $ on
$ {E}_{s}^{1} := {}^{p}\varphi _{t_{s},1}L[1] $ is semisimple so that
$ {\rho }_{s}^{1} = {\mu }_{s}^{1} + {\nu }_{s}^{1} $,
$ {\rho }_{s}^{\prime 1} = {\nu }_{s}^{1} $,
where
$ {}^{p}\varphi _{t_{s},1} $ is the unipotent monodromy part of
$ {}^{p}\varphi _{t_{s}} $.

\bigskip
\noindent
{\bf 3.7.}
{\it Remark.} In the case
$ X = {\Bbb A}^{n}, S = {\Bbb A}^{1} $,
we can calculate the Gauss-Manin system and the Brieskorn-Deligne
lattice in the following way provided that the polynomial
$ f $ is not very complicated :

\noindent
(i) Calculate
$ df\wedge d\Omega ^{n-2} $ to get a basis of
$ {G}_{f}^{(0)} = \Omega ^{n}/df\wedge d\Omega ^{n-2 }  $ over
$ {\Bbb C} $.

\noindent
(ii) Calculate the action of
$ t $ on the basis to get generators of
$ {G}_{f}^{(0)} $ over
$ {\Bbb C}[t] $.

\noindent
(iii) Calculate
$ {\partial }_{t}^{-1} $ (using
$ {\partial }_{t}^{-1}(\omega ) = df\wedge \omega ' $ for
$ d\omega ' = \omega ) $ to get the differential equations and 
determine
$ {G}_{f}^{(0),>0} $.

\noindent
Here
$ \Omega ^{p} = \Gamma (X, {\Omega }_{X}^{p}) $ as in the introduction.
The argument is similar to [7] except that our case is algebraic and 
global.
For example, we use
$ {\partial }_{t}^{-1} $ rather than
$ \partial _{t} $,
and
$ (\prod _{c} (t-c)^{k_{c}})\partial _{t} $ is used at the last stage only 
for a direct factor or a subquotient of
$ G_{f} $.

\bigskip
\noindent
{\bf 3.8.}
{\it Examples.} (i)  $ f = y^{2} + x^{3} - 3x $.
This is the simplest example such that
$ \bar{L} \ne  0 $ in the notation of (0.8).
We have
$ h = g = 1 $,
$ m = 2 $,
$ r = \nu  = 0 $,
$ \mu  = 2 $,
$ {\mu }_{\pm 2}^{1} = 1 $,
$ R'_{\pm 2} = 2, R_{\pm 2} = 1 $,
$ N'_{\pm 2} = N_{\pm 2} = 0 $,
and get the Gauss hypergeometric differential equation.  
The calculation seems much easier than the conventional one using the 
discriminant and the integration by parts as explained in the 
introduction of [21].
The indicial equation [8] is compatible with
$ {G}_{f}^{(0)} = {G}_{f}^{>0} $.

(ii)  $ f = x^{4}y^{2} + 2x^{2}y + xy^{2} $.
This is an example such that
$ h'' \,(= \dim H^{1}(S^{\an}, j_{*}\bar{L})) \ne  0 $.
We have
$ m = 7 $,
$ g = 1 $,
$ h = 5 $,
$ r = 1 $,
$ h' = h'' = 2 $,
$ \mu  = \mu _{0} = 4 $,
$ \nu  = \nu _{-1} = 3 $,  (where
$ \nu  $ can be calculated as in [15] combined with [19]), and
$ R'_{0} = 7 $,
$ R_{0} = 5 $,
$ N'_{0} = 0 $,
$ N_{0} = 2 $,
$ R'_{-1} = R_{-1} = 4 $,
$ N'_{-1} = N_{-1} = 0 $.

(iii)  $ f = x^{2}y^{2} + 2xy + x $.
This is an example such that
$ \nu  \ne  0 $ and
$ \bigcap _{i} {\partial }_{t}^{-i}{G}_{f}^{(0)} = 0 $.
We have
$ h = m = 2 $,
$ h' = r = 1 $,
$ g = 0 $,
$ \mu  = {\mu }_{0}^{1} = 1 $,
$ \nu  = {\nu }_{-1}^{-1} = 1 $,
and
$ {\Cal G}_{f} = 
{\Cal D}_{S}/{\Cal D}_{S}t\partial _{t}((t+1)\partial _{t} + 1/2) $.
In this case
$ {G}_{f}^{(0)} \rightarrow  G_{f}/_{\tor}G_{f} $ is bijective.

Note that the condition
$ \mu  > 0 $ does not necessarily imply that
$ \bigcap _{i\in {\Bbb N}} {\partial }_{t}^{-i}{G}_{f}^{(0)} = 0 $ (e.g.
 $ f = x^{2}y^{2} + x^{2} + 2x $ or
$ f = x^{6}y^{3} + x^{2}y^{2} + 2xy) $.
This property holds for the local analytic Brieskorn lattice, and 
hence it is true if we take a completion of
$ {\Cal O} $ by some topology associated with
$ x \in  \Sing f $.
But it is not true without taking it.

\bigskip\bigskip
\centerline{{\bf References}}

\bigskip

\item{[1]}
E.~Artal Bartolo, P.~Cassou-Nogu\`es and A.~Dimca, Sur la topologie des 
polyn\^omes complexes, in Proceedings Oberwolfach Singularities 
Conference 1996, Brieskorn Festband, Progress in Math., vol.~162, 
Birkh\"auser, Basel, 1998, pp.~317--343.

\item{[2]}
G.~Bailly-Maitre, Sur le syst\`eme local de Gauss-Manin d'un polyn\^ome 
de deux variables, preprint, Bordeaux University, 1998.

\item{[3]}
A.~Beilinson, J.~Bernstein and P.~Deligne, Faisceaux Pervers, 
Ast\'erisque, vol.~100, Soc.~Math.~France, Paris, 1982.

\item{[4]}
A.~Borel et al., Algebraic $ {\Cal D} $-modules, Perspectives in Math.~2,
Academic  Press, 1987.

\item{[5]}
L.~Boutet de Monvel, ${\Cal D}$-modules holon\^omes r\'eguliers en une 
variable, in Math\'ematique et Physiques, Progress in Math., 
Birkh\"auser, vol.~37, (1983), pp.~281--288.

\item{[6]}
J.~Brian\c con and Ph.~Maisonobe, Id\'eaux de germes d'op\'erateurs 
diff\'erentiels \`a une variable, Enseign.~Math., 30 (1984), 7--36.

\item{[7]}
E.~Brieskorn, Die Monodromie der isolierten Singularit\"aten von 
Hyperfl\"achen, Manuscripta Math., 2 (1970), 103--161.

\item{[8]}
E.~Coddington, N.~Levinson, Theory of Ordinary Differential Equations, 
McGraw-Hill, 1955.

\item{[9]}
P.~Deligne, Equations Diff\'erentielles \`a Points Singuliers R\'eguliers, 
Lect. Notes in Math. vol.~163, Springer, Berlin, 1970.

\item{[10]}
\SameAuthor, Le formalisme des cycles \'evanescents, in SGA7 XIII and XIV, 
Lect. Notes in Math. vol.~340, Springer, Berlin, 1973, pp. 82--115 and 
116--164.

\item{[11]}
\SameAuthor, Th\'eorie de Hodge I, Actes Congr\`es Intern. Math., 1970, 
vol. 1, 425-430 : II,  Publ. Math. IHES, 40 (1971), 5--57; III ibid., 44 
(1974), 5--77.

\item{[12]}
A.~Dimca, A sheaf-theoretic view at the topology of polynomials, 
Proc. Conf. on Algebraic Geometry in Constantza, 1996, 
An. St. Univ. Ovidius Constantza, 5 (1997), 17--22.

\item{[13]}
\SameAuthor, Invariant cycles for complex polynomials,
 Rev. Roumaine Math. Pures Appl. 43 (1998), 113--120.

\item{[14]}
\SameAuthor, Monodromy at infinity for polynomials in two variables, 
Journal of Algebraic Geometry 7 (1998) 771--779.

\item{[15]}
A.~Durfee, Five definitions of critical points at infinity, in Proceedings 
Oberwolfach Singularities Conference 1996, Brieskorn Festband, 
Progress in Math., vol. 162, Birkh\"auser, Basel, 1998, pp. 345--360.

\item{[16]}
L.~Gavrilov, Petrov modules and zeros of abelian integrals, 
Bull. Sci. Math., 122 (1998), 571--584.

\item{[17]}
S.~Kaliman, Two remarks on polynomials in two variables, 
Pacific J. Math. 154 (1992), 285--295.

\item{[18]}
M.~Kashiwara, Vanishing cycle sheaves and holonomic systems of
differential equations, Lect. Notes in Math., vol. 1016, Springer, Berlin,
 1983, pp. 136--142 .

\item{[19]}
A.~Kouchinirenko, Poly\`edres de Newton et nombres de Milnor, 
Inv. Math. 32 (1976), 1--31.

\item{[20]}
B.~Malgrange, Polyn\^ome de Bernstein-Sato et cohomologie \'evanescente, 
Ast\'erisque, 101--102 (1983), 243--267.

\item{[21]}
F.~Pham, Singularit\'es des Syst\`emes Diff\'erentiels de Gauss-Manin, 
Progress in Math. vol. 2, Birkh\"auser, Basel, 1979.

\item{[22]}
\SameAuthor, Structures de Hodge mixtes associ\'ees \`a un germe de fonction 
\`a point critique isol\'e, Ast\'erisque, 101--102 (1983), 268--285.

\item{[23]}
C.~Sabbah, Hypergeometric periods for a tame polynomial,
C. R. Acad. Sci. Paris, 328 S\'erie I (1999) 603--608.

\item{[24]}
M.~Saito, Modules de Hodge polarisables, Publ. RIMS, Kyoto Univ., 24 
(1988), 849--995.

\item{[25]}
\SameAuthor, Mixed Hodge Modules, Publ. RIMS Kyoto Univ. 26 (1990), 
221--333.

\item{[26]}
\SameAuthor, On the structure of Brieskorn lattice, Ann. Institut Fourier 39 
(1989), 27--72.

\item{[27]}
\SameAuthor, Hodge filtrations on Gauss-Manin systems I, 
J. Fac. Sci. Univ. Tokyo, Sect. I A 30 (1984), 489--498; 
II, Proc. Japan Acad., Ser. A 59 
(1983), 37--40.

\item{[28]}
\SameAuthor, Period mapping via Brieskorn modules, Bull. Soc. Math. France, 
119, (1991), 141--171.

\item{[29]}
J.~Scherk and J.~Steenbrink On the mixed Hodge structure on the 
cohomology of the Milnor fiber, Math. Ann. 271 (1985) 641--665.

\item{[30]}
M.~Sebastiani, Preuve d'une conjecture de Brieskorn, Manuscripta Math., 
2 (1970), 301--308.

\item{[31]}
D.~Siersma and M.~Tib\u ar, Singularity at infinity and their vanishing 
cycles, Duke Math. J. 80 (1995), 771--783.

\item{[32]}
J.~Steenbrink, Mixed Hodge structure on the vanishing cohomology, 
in Real and Complex Singularities (Proc. Nordic Summer School, Oslo, 
1976) Alphen a/d Rijn: Sijthoff \& Noordhoff 1977, pp. 525--563.

\item{[33]}
M.~Suzuki, Propri\'et\'es topologiques des polyn\^ome de deux variables 
complexes et automorphismes alg\'ebriques de l'espace
$ C^{2} $, J. Math. Soc. Japan 26 (1974), 141--157.

\item{[34]}
M.~Tib\u ar, Topology at infinity of polynomial mappings and Thom 
regularity condition, Compos. Math. 111 (1998), 89--109.

\item{[35]}
A.~Varchenko, The asymptotics of holomorphic forms determine a mixed 
Hodge structure, Soviet Math. Dokl., 22 (1980), 772--775.

\item{[36]}
\SameAuthor,
Asymptotic mixed Hodge structure in vanishing cohomology, 
Math. USSR Izvestija 18 (1982), 469--512.

\item{[37]}
S.~Zucker, Hodge theory with degenerating coefficients, 
$ L_{2} $ cohomology in the Poincar\'e metric, Ann. Math. 109
(1979), 415--476.
\bigskip

\noindent
Alexandru Dimca

\noindent
Math\'ematiques pures, Universit\'e Bordeaux I

\noindent
33405 Talence Cedex, FRANCE

\noindent
e-mail: dimca\@math.u-bordeaux.fr

\bigskip

\noindent
Morihiko Saito

\noindent
RIMS Kyoto University

\noindent
Kyoto 606--8502 JAPAN

\noindent
e-mail: msaito\@kurims.kyoto-u.ac.jp

\bye